\documentstyle[12pt,amssymb]{article}

\newcommand{\sector}[1]{\section{#1}\setcounter{thm}{0}}

\newtheorem{thm}{Theorem}[section]
\newcommand{\theorem}{\begin{thm}}
\newcommand{\theoremEnd}{\end{thm}}

\newtheorem{pro}[thm]{Proposition}
\newcommand{\proposition}{\begin{pro}}
\newcommand{\propositionEnd}{\end{pro}}

\newtheorem{cor}[thm]{Corollary}
\newcommand{\corollary}{\begin{cor}}
\newcommand{\corollaryEnd}{\end{cor}}

\newtheorem{lem}[thm]{Lemma}
\newcommand{\lemma}{\begin{lem}}
\newcommand{\lemmaEnd}{\end{lem}}

\newtheorem{cnj}[thm]{Conjecture}
\newcommand{\conjecture}{\begin{cnj}}
\newcommand{\conjectureEnd}{\end{cnj}}

\newtheorem{dft}[thm]{Definition}
\newcommand{\definition}{\begin{dft}}
\newcommand{\definitionEnd}{\end{dft}}

\newcommand{\proof} {\noindent {\em Proof. }}
\newcommand{\proofSketch} {\noindent {\em Proof Sketch. }}
\newcommand{\qed} {\hfill$\Box$\bigskip}
\newcommand{\qedNoSkip} {\hfill$\Box$}
\newcommand{\qedInEquation} {\eqno{\Box}$$}
\newcommand{\proofOf}[1] {\noindent {\em Proof of #1.}}

\newcommand{\union}{\bigcup}
\newcommand{\intt}{\bigcap}
\newcommand{\eps}{\varepsilon}
\newcommand{\reals}{{\Bbb R}}

\newcommand{\eqdef}{\;\;{:=}\;}

\newcommand{\Hom}{\mbox{Hom}}
\newcommand{\End}{\mbox{End}}
\newcommand{\Ker}{\mbox{Ker}}
\newcommand{\Tr}{\mbox{Tr}}

\title{Circle-valued Morse theory, Reidemeister torsion, and Seiberg-Witten
invariants of 3-manifolds} 
\author{Michael Hutchings \& Yi-Jen Lee}
\date{First version: Sep.\ 6, 1996\\ This version: Dec.\ 1, 1996}

\begin{document}

\maketitle

\begin{abstract}
Let $X$ be a compact oriented Riemannian manifold and let $\phi:X\to S^1$
be a circle-valued Morse function.  Under some mild assumptions on $\phi$,
we prove a formula relating: 
\begin{description}
\item{(a)}
the number of closed orbits of the gradient flow of $\phi$ of any given degree;
\item{(b)}
the torsion of a ``Morse complex'', which counts gradient flow lines
between critical points of $\phi$; and 
\item{(c)} 
a kind of Reidemeister torsion of $X$ determined by the homotopy class of
$\phi$.
\end{description}
When $\dim(X)=3$ and $b_1(X)>0$, we state a conjecture analogous to
Taubes's ``SW=Gromov'' theorem, and we use it to deduce (for closed
manifolds, modulo signs) the Meng-Taubes relation between the
Seiberg-Witten invariants and the ``Milnor torsion'' of $X$.
\end{abstract}

\tableofcontents

\sector{Introduction and statement of results}

It has long been known that one can obtain information about the homology
of a manifold from the structure of the critical points of a Morse function
defined on it.  A sharp statement of this relationship is that there
is an isomorphism between the homology of the manifold and the homology of
a ``Morse complex'' whose chains are critical points and whose differential
counts gradient flow lines between critical points.

Novikov \cite{novikov} generalized this relationship to circle-valued (and
other multiply-valued) Morse functions.  One can still define a Morse
complex in terms of gradient flow lines between critical points, which is
now a module over the ring of integer Laurent series in one variable.
Novikov obtained bounds on the betti numbers of the manifold in terms of
algebraic invariants of the homology of this complex.

The novelty of this paper is that we also consider the {\em closed
orbits\/} of the gradient flow of a circle-valued Morse function.  These
turn out to be related not to homology, but rather to Reidemeister torsion.
(For nonsingular flows, relations between closed orbits and Reidemeister
torsion have been investigated by D.\ Fried \cite{fried:homological}
\cite{fried:lefschetz}.)

\subsection{Counting closed orbits}

Let $X^n$ be a compact oriented Riemannian manifold and let $\phi:X\to
S^1={\Bbb R}/{\Bbb Z}$ be a generic Morse function.  (See
\S\ref{sec:morseHomology} for the definition of ``generic''.)  Assume $0\in
{\Bbb R}/{\Bbb Z}$ is a regular value of $\phi$, and let
$\Sigma=\phi^{-1}(0)$.  If $p\in\Sigma$, we can flow upwards from $p$ along
the gradient vector field of $\phi$.  As $\phi$ goes once around $S^1$, we
will return to $\Sigma$, if we do not get sucked into a critical point
first.  Let $f(p)$ denote this point of $\Sigma$, if it exists.  Thus $f$
is a function from a subset of $\Sigma$ to a subset of $\Sigma$.  Also
$f^k$ is a function defined on a smaller subset of $\Sigma$.  Let
$\mbox{Fix}(f^k)$ denote the signed number of fixed points of $f^k$.  (A
fixed point of a function corresponds to an intersection point of the graph
with the diagonal, and the sign of this intersection number determines the
sign of the fixed point.)

The Morse complex
$$M^0\stackrel{d}{\longrightarrow}
M^1\stackrel{d}{\longrightarrow}\cdots\stackrel{d}{\longrightarrow}M^n$$
is defined as follows.  Let $L_{{\Bbb Z}}$ be the ring of Laurent series in
one variable $t$ with integer coefficients, i.e.\ formal sums
$\sum_{k=k_0}^\infty a_kt^k$ with $a_k\in{\Bbb Z}$.  Let $M^i$ be the free
$L_{{\Bbb Z}}$-module generated by $\mbox{Crit}^i$, the set of critical
points of index $i$.  If $x\in\mbox{Crit}^i$, define
$$dx\eqdef\sum_{y\in\mbox{\scriptsize Crit}^{i+1}}\langle dx,y\rangle y,$$
where $\langle dx,y\rangle$ is a Taylor series whose $n^{th}$ coefficient
is the signed number of gradient flow lines from $x$ to $y$ that cross
$\Sigma$ $n$ times.  (The sign conventions, and other details, are
explained in \S\ref{sec:morseHomology}.)

Let us recall the definition of Reidemeister torsion.  Suppose
$C^0\stackrel{d}{\longrightarrow} 
C^1\stackrel{d}{\longrightarrow}\cdots\stackrel{d}{\longrightarrow} C^m$ is
an acyclic complex of finite dimensional vector spaces over a field $F$,
and suppose that each 
vector space $C^i$ has a volume form chosen on it.  Choose
$\omega_i\in\wedge^*C^i$, $i=0,\ldots,m-1$, so that
$d\omega_{i-1}\wedge\omega_{i}\in\wedge^{\mbox{\scriptsize
top}}C^i$.  Then the Reidemeister torsion $\tau(C)$ is defined to be
\begin{equation}
\label{eqn:tauDef}
\tau(C)\eqdef\prod_{i=0}^m\mbox{vol}(d\omega_{i-1}\wedge\omega_i)^{(-1)^i}
\end{equation}
(where we interpret $d\omega_{-1}=1$).  One can check, using the fact that
$d^2=0$, that $\tau(C)$ does not depend on the choice of $\omega_i$'s.  If
the differential has degree -1 instead of 1, the definition is analogous,
but we adopt the following sign
convention: if we have a complex $C_m\stackrel{\partial}{\longrightarrow}
C_{m-1}\stackrel{\partial}{\longrightarrow} \cdots
\stackrel{\partial}{\longrightarrow} C_0$, then we choose
$\omega_i\in\wedge^*C_i$, $i=1,\ldots,m$, so that 
$\partial\omega_{i+1}\wedge\omega_{i}\in\wedge^{\mbox{\scriptsize
top}}C_i$ and define
$$
\tau(C) \eqdef
\prod_{i=0}^m\mbox{vol}(\partial\omega_{i+1}\wedge\omega_i)^{(-1)^{m-i}}. 
$$

There is a natural symmetric bilinear form $\langle\;,\;\rangle$ on the
Morse complex, in which the critical points are orthonormal.  This defines
a volume form on the $L_{{\Bbb Q}}$-vector space $M^i\otimes L_{{\Bbb Q}}$.
(Here $L_{{\Bbb Q}}$ is the field of rational Laurent series.)  If the
complex $M^*\otimes L_{{\Bbb Q}}$ is acyclic, then the Reidemeister torsion of
$M^*\otimes L_{{\Bbb Q}}$ is defined, up to sign.  We denote this
simply by $\tau(M)$.  (One can define this in the non-acyclic case as well,
but we will not need to.)

Let $\tilde{X}$ be the infinite cyclic cover of $X$ induced by $\phi$,
i.e.\ the fiber product of $X$ and ${\Bbb R}$ over $S^1$: 
$$
\begin{array}{ccc}
\tilde{X}&\longrightarrow&{\Bbb R}\\
\big\downarrow&&\big\downarrow\\
X&\stackrel{\phi}{\longrightarrow}&{\Bbb R}/{\Bbb Z}.
\end{array}
$$
We think of $\tilde{X}$ as a subset of $X\times{\Bbb R}$.  There is a
covering transformation of $\tilde{X}$ which shifts points ``down'',
i.e. which sends $(p,\lambda)\mapsto(p,\lambda-1)$ for $p\in X$ and
$\lambda\in\reals$.  Let $A:H_*(\tilde{X};{\Bbb Q})\to
H_*(\tilde{X};{\Bbb Q})$ be the map in rational homology induced by this
covering transformation.

We can now state the formula for counting closed orbits.

\theorem
\label{thm:main}
Assume $H^*(M^*\otimes L_{{\Bbb Q}})=0$.  On $\partial X$, assume that
$\mbox{grad}(\phi)$ is parallel to the boundary and has no zeroes there.
Then
$$\sum_{k=1}^\infty t^k\mbox{Fix}(f^k)-(-1)^n t\frac{d}{dt}\log\tau(M) =
\Tr(tA(1-tA)^{-1})+m,$$
where $m\in{\Bbb Z}$.
\theoremEnd

\paragraph{Notes.}
Here $\Tr$ denotes the graded trace.
We will see in Corollary~\ref{cor:surjective} that $H_*(\tilde{X};{\Bbb Q})$
is finite dimensional (thanks to our assumption that
$H^*(M^*\otimes L_{{\Bbb Q}})=0$), so this trace is well defined.  Also
$\frac{d}{dt}\log\tau$ means $\tau^{-1}\frac{d}{dt}\tau$, which is well
defined even though $\tau$ has a sign ambiguity.

\paragraph{Example.} If there are no critical points, then $\tau(M)=1$, $f$
is a diffeomorphism of $\Sigma$, and our theorem reduces to the Lefschetz
fixed point formula for $f$.  Thus we can think of the theorem as Lefschetz
fixed point formula for certain partially defined functions, in which the
R-torsion of the Morse complex appears as a correction term.

\subsection{A refinement}

This theorem is the logarithmic derivative of a slightly sharper
formula, which we will also prove.  To state it we need two more
definitions.  First, define a ``zeta function''
\begin{equation}
\label{eqn:zetaDef}
\zeta(f)  \eqdef  \exp\left(\sum_{k\ge 1}\mbox{Fix}(f^k)\frac{t^k}{k}\right).
\end{equation}
Note that
$$\zeta(f) = \prod_{\gamma\in{\cal O}}(1-t^{k(\gamma)})^{-\eps(\gamma)},
$$
where ${\cal O}$ is the set of irreducible, connected closed orbits, and for
$\gamma\in{\cal O}$, $k(\gamma)$ is the degree of $\phi:\gamma\to S^1$ and
$\eps(\gamma)$ is the
sign of each of the corresponding $k$ fixed points of $f^k$. 

Second, let $C_*(\tilde{X})$ be the chain complex associated to a cell
decomposition of $\tilde{X}$ lifted from a cell decomposition of $X$. This
is a ${\Bbb Z}[t,t^{-1}]$ module, where $t$ acts via the upward covering
transformation.  Let $Q({\Bbb Z}[t,t^{-1}])$ be the field of fractions of
${\Bbb Z}[t,t^{-1}]$.  Our assumption on the acyclicity of the rational
Morse complex implies that $C_*(\tilde{X})\bigotimes_{{\Bbb Z}[t,t^{-1}]}
Q({\Bbb Z}[t,t^{-1}])$ is acyclic.  (This follows from
Corollary~\ref{cor:isomorphism}.) Furthermore this complex has a volume
form, well defined up to multiplication by $\pm t^k$, consisting of a wedge
product of cells, one for each cell in $X$.  Thus the torsion of this
complex, which we denote by $\tau(X,\phi)$, is a well defined element of
$Q({\Bbb Z}[t,t^{-1}])/\pm t^k$.  (Of course the homology of this complex is
isomorphic to $H_*(\tilde{X})$, viewed as a ${\Bbb Z}[t,t^{-1}]$ module.  We
will see in Lemma~\ref{lem:torsion} that this homology alone, and not the
choice of cell decomposition, determines the torsion.)

\paragraph{Example.} Let $K\subset S^3$ be a knot, and let $X$ be the
3-manifold 
obtained from $(0,1)$ surgery on $K$.  Let $\phi:X\to S^1$ be a Morse
function whose homotopy class in $H^1(X;{\Bbb Z})$ is the Alexander dual of
$K$, i.e. which sends a loop $\gamma$ in $X$ to the linking number of
$\gamma$ with $K$.  Then
$$\tau(X,\phi)=\frac{\Delta_K(t)}{(1-t)^2},$$
where $\Delta_K(t)$ is the Alexander polynomial of $K$.

\theorem
\label{thm:refinement}
Under the assumptions of Theorem~\ref{thm:main}, we have
$$\zeta(f)^{(-1)^{n-1}}\tau(M)=\tau(X,\phi),$$
up to multiplication by $\pm t^k$.
\theoremEnd

\paragraph{Remark.} The right side of this formula clearly depends only on
the homotopy class of $\phi$ in $H^1(X;{\Bbb Z})$.  It is not too hard to
show directly that the left hand side is also invariant, but this requires
our assumption that the rational Morse complex is acyclic.  To see this,
note that the zeta function is determined by the intersection number of the
graph of $f^k$ and the diagonal in $\Sigma\times\Sigma$, for
$k=1,2,\ldots$.  If we deform $\phi$, this intersection number will change
exactly when the boundary of the graph crosses the diagonal.  But we will
see that the boundary of the graph consists of points of the form $(p,q)$
where $p$ is in the descending manifold of a critical point $x$, and $q$ is
in the ascending manifold of $x$.  When the ascending and descending
manifolds cross each other this way, the effect is to replace $dx$ by
$(1-t^k)^{\pm 1}dx$ and $d^*x$ by $(1-t^k)^{\mp1}d^*x$.  In the acyclic
case, this multiplies or divides $\tau(M)$ by $(1-t^k)$, which exactly
cancels the change in the zeta function.  But in the nonacyclic case, if
say $dx=d^*x=0$, then $\tau(M)$ does not change, even though $\zeta(f)$
changes.

Of course this formula might generalize to the nonacyclic case if a
suitable correction term is included.  However several of the steps of our
proof do not appear to have natural generalizations to the nonacyclic case.

\paragraph{Example.} Suppose $X=S^1$ and $\phi:S^1\to S^1$ has degree $k$.
Then $\tilde{X}$ has $k$ components, and the upward covering transformation
permutes them in a $k$-cycle, so we find that $\tau(X,\phi)=(1-t^k)^{-1}$.  If
$\phi$ has no critical points then $\zeta(f)=(1-t^k)^{-1}$ and $\tau(M)=1$.  If
there are $c>0$ critical points then $\zeta(f)=1$, while the differential
$d$ of 
the Morse complex has the form
$$
\left(
\begin{array}{ccccc}
t^{a_1} & -t^{b_1} & 0 & \cdots & 0 \\
0 & t^{a_2} & -t^{b_2} & \cdots & 0 \\
\vdots & & \ddots & & \vdots \\
-t^{b_c} & 0 & \cdots & 0 & t^{a_c}
\end{array}
\right)
$$
(for a certain choice of orientations).  Since $\phi$ has degree $k$, we have
$\sum a_i-\sum b_i=k$, so up to signs and powers of $t$,
$\tau(M)=(1-t^k)^{-1}$.

\subsection{Relation to Seiberg-Witten theory}

Suppose now that $X$ is a closed oriented 3-manifold with $b^1(X)>0$.  Let
${\cal S}$ denote the set of $Spin^c$ structures on $X$.  A $Spin^c$
structure is a $U(2)$ bundle $S\to X$ such that the fiber over each point
is an irreducible Clifford module over the tangent space.  The set
${\cal S}$ is an $H^2(X;{\Bbb Z})$-torsor;  $\alpha\in H^2(X;{\Bbb Z})$ sends
$S$ to $S\otimes{\cal L}$, where ${\cal L}$ is the complex line bundle with
$c_1({\cal L})=\alpha$.  Suppose an orientation is chosen on the rational
homology of $X$.  Then the Seiberg-Witten invariant
$$SW:{\cal S}\to{\Bbb Z}$$
is defined.  (See e.g.\ \cite{kronheimer-mrowka}, \cite{morgan-szabo-taubes},
\cite{meng-taubes}, or \cite{witten:monopoles}, \cite{morgan} for the four
dimensional case.)  (When $b^1(X)=1$, the definition of $SW$ requires some
care; see \S\ref{sec:conjecture}.)

Taubes \cite{taubes:SWGromov} has shown that for a symplectic 4-manifold
with a metric compatible with the symplectic form, the Seiberg-Witten
invariant is equal to the Gromov invariant (see Taubes
\cite{taubes:counting} for a precise definition), which counts
pseudoholomorphic curves.  He has also made some progress on generalizing
this result to nonsymplectic 4-manifolds \cite{taubes:lectures}.  Here one
considers 2-forms which are symplectic except on a set of circles, and
pseudoholomorphic curves bounded by this set of circles.

In three dimensions, a possible analogue of a symplectic form is a harmonic
1-form, and the analogue of pseudoholomorphic curves is flow lines of the
dual vector field.  Actually we do not need to assume that the 1-form is
harmonic, but only that it has no index 0 or 3 critical points.  In
\S\ref{sec:conjecture} we define an analogue of the Gromov invariant out of
such a 1-form.  This turns out to be a map
$$I:{\cal S}\to{\Bbb Z}.$$
It is similar to the left hand side of Theorem~\ref{thm:refinement}, but it
is sharper, because it keeps track of the relative homology classes of
closed orbits and gradient flow lines (and not just the intersection
numbers with $\Sigma$).
A priori $I$ depends on the choice of an integral cohomology class, but we
conjecture (Conjecture~\ref{cnj:main}) that it does not, and in fact
$$SW=\pm I.$$
(The basic idea of this was suggested to us by Taubes.)

In \S\ref{sec:meng-taubes}, we apply Theorem~\ref{thm:refinement} to prove:

\theorem[assuming Conjecture~\ref{cnj:main}]
\label{thm:application}
Let $X$ be a closed oriented 3-manifold with $b^1>0$ and $0\neq\alpha\in
H^1(X;{\Bbb Z})$.  Then
$$
\sum_{S\in{\cal S}}SW(S)t^{\alpha(c_1(\det S))/2} =
\left\{ \begin{array}{cl}
\tau(X,\phi) & \mbox{if $M^*\otimes L_{{\Bbb Q}}$ acyclic}\\ 
0 & \mbox{otherwise}
\end{array}
\right.
$$
modulo multiplication by $\pm t^k$, where $\phi:X\to S^1$ is in the
homotopy class determined by $\alpha$.
\theoremEnd

When $b^1(X)>1$, the $t^k$ ambiguity in the right hand side of this theorem
can be resolved by applying the ``charge conjugation invariance'' of the
Seiberg-Witten equations, which tells us that the left hand side is
invariant under $t\mapsto t^{-1}$.  (See Witten~\cite{witten:monopoles} or
Morgan~\cite{morgan} for the 4-dimensional case.  The 3 dimensional case is
also easy to deduce from Conjecture~\ref{cnj:main}, by replacing $\eta$ in
\S\ref{sec:conjecture} with $-\eta$.)

Note that even if we apply this theorem to every $\alpha\in
H^1(X;{\Bbb Z})$, we cannot recover all of the Seiberg-Witten invariants of
$X$, because the theorem does not distinguish between $Spin^c$ structures
that differ by the action of a torsion element of $H^2(X;{\Bbb Z})$.
However we can recover the theorem of Meng and Taubes \cite{meng-taubes}
relating $SW$ to ``Milnor torsion'' (modulo signs, in the case of closed
manifolds). We explain this in \S\ref{sec:meng-taubes}.  (Actually our
methods might also be applicable 
to the Meng-Taubes formula for manifolds with boundary,
because the boundary conditions in Theorem~\ref{thm:main} and in
\cite{meng-taubes} are similar.)

\paragraph{Acknowledgments}

We are grateful to Cliff Taubes for
pointing us in this direction, for helpful conversations, and for his
advice and support in general.  The first
author would like to thank Raoul Bott, Ken Fan and Ravi Vakil for
additional helpful conversations, and Jim Bryan for teaching him Morse
theory long ago.
The second author would like to thank Peter Kronheimer for additional
helpful conversations.

\sector{Preliminaries}

\subsection{Notation and conventions}

The symbol $\cdot$ always denotes intersection number.  We take
$\alpha\cdot\beta$ to be zero if $\alpha$ and $\beta$ do not have
complementary dimension.

The differential in the Morse complex is $d$, its adjoint with respect to
the natural inner product $\langle\;,\;\rangle$ is $d^*$, and $\Delta\eqdef
dd^*+d^*d$.

The symbol `Tr' always denotes graded trace.

If $R$ is an integral domain, $Q(R)$ is its field of fractions, and $L_R$
is the ring of Laurent series with coefficients in $R$, i.e.\ functions
${\Bbb Z}\to R$ supported away from $-\infty$.

Regarding orientations: if $Y$ is an oriented manifold, $\psi$ is a Morse
function on $Y$, and
$\lambda\in{\Bbb R}$ is a regular value of $\psi$, we orient the level set
$W=\psi^{-1}(\lambda)$ by declaring
$$TY|W={\Bbb R}\cdot\mbox{grad}(\psi)\oplus TW$$
to be an isomorphism of oriented vector bundles.  If $Y$ is an oriented
manifold with boundary, we orient the boundary via the convention
$$TY|\partial Y={\Bbb R}\cdot\nu\oplus T(\partial Y),$$
where $\nu$ points outwards.

\subsection{Morse homology}
\label{sec:morseHomology}

We will now present a direct approach to the relation between
Morse homology and ordinary homology.  We will then deduce some simple lemmas
concerning circle-valued Morse homology.  

The approach here is very useful for understanding Theorem~\ref{thm:main},
but there is one technical difficulty, which is that we have to compute
ordinary homology using special singular chains.  These chains must have
well defined intersection numbers with the descending manifolds of the
critical points, and the set of such chains must be preserved by the
gradient flow.  For example, we can use smooth chains with conical
singularities, in which the smooth set and the singular sets intersect the
descending manifolds of the critical points transversely.  (See Laudenbach
\cite{laudenbach} for a proof that these are preserved by gradient flow.)
One can use standard techniques to show that an arbitrary chain may be
approximated by these special chains, so that the homology of the complex
of special chains is isomorphic to ordinary homology.  We will not go into
further details about this.

Let $Y^n$ be an oriented Riemannian manifold and let $\psi:Y\to(-\infty,0]$
be a generic Morse function with $\partial Y=\psi^{-1}(0)$. (``Generic''
means that the gradient flow is Morse-Smale, i.e.\ the ascending and
descending manifolds of different critical points intersect transversely.)
Let $\mbox{Crit}$ denote the set of critical points.  For
$x\in\mbox{Crit}$, we define the ascending manifold of $x$,
${\cal A}(x)$, to be the closure of the set of all $p\in Y$ such that
downward gradient flow from $p$ converges to $x$.  We define the descending
manifold ${\cal D}(x)$ analogously, using upward gradient flow.  Choose
orientations on ${\cal A}(x)$ and ${\cal D}(x)$ such that the
intersection number ${\cal A}(x)\cdot{\cal D}(x)$ is $+1$ in $Y$.

We define the Morse complex as follows.  Let $M^i$ be the free
${\Bbb Z}$-module generated by $\mbox{Crit}^i$, the set of critical points
of index $i$.  Define $d:M^i\to M^{i+1}$ by
$$dx\eqdef\sum_{y\in\mbox{\scriptsize Crit}^{i+1}}\langle dx,y\rangle y,$$
where $\langle dx,y\rangle$ is the signed number of gradient flow lines
from $x$ to $y$.  A flow line counts with positive sign when the
intersection number $({\cal D}(x)\cap\psi^{-1}(\lambda))
\cdot ({\cal A}(y)\cap\psi^{-1}(\lambda))$ is $+1$ in $\psi^{-1}(\lambda)$.

If $\alpha\in C_*(Y)$ is a generic chain, define ${\cal F}(\alpha)$ to
be the closure of the union, over all $s\in[0,\infty)$, of the time $s$
upward gradient flow applied to $\alpha$, oriented so that the
orientations on $\partial{\cal F}(\alpha)$ and $\alpha$ disagree.  We need
just one geometric observation:

\lemma
\label{lem:flowGeometry}
If $\alpha$ is a generic chain, then
$$\partial{\cal F}(\alpha)={\cal F}(\alpha)\cap\partial Y -\alpha
-{\cal F}(\partial \alpha)+\sum_{x\in\mbox{\scriptsize
Crit}}({\cal D}(x)\cdot\alpha){\cal A}(x).$$
\lemmaEnd

\proofSketch
This is straightforward, using the Morse lemma to model the behavior near
the critical points.  (It is a little easier to first assume all the
critical points are at the same height, and then use induction.)  The idea
is illustrated in Figure~1.
\qedNoSkip

\lemma
\label{lem:AscChainMap}
If $x$ is a critical point then
$$
\partial{\cal A}(x) = {\cal A}(x)\cap\partial Y+{\cal A}(dx).
$$
(Here we are extending ${\cal A}$ to a map from $M^*$ to
$C_{n-*}(Y)$.)
\lemmaEnd

\proof
Let $s\in{\Bbb R}$ be slightly larger than $\psi(x)$, and let
 $\alpha={\cal A}(x)\intt\psi^{-1}(s)$.  Apply
 Lemma~\ref{lem:flowGeometry} to obtain
$$\partial({\cal A}(x)\cap\psi^{-1}[s,0])={\cal A}(x)\cap\partial
Y-\alpha+{\cal A}(dx).$$
Since $\partial({\cal A}(x)\cap\psi^{-1}(-\infty,s])=\alpha$,
we are done.
\qed

It follows from this lemma that
$${\cal A}:M^*\to C_{n-*}(Y,\partial Y)$$
is a chain map.  It then follows from $\partial^2=0$
that
$$d^2=0,$$
since the ascending manifolds of different critical points are disjoint.

\proposition
\label{pro:morseHomology}
The chain map ${\cal A}$ induces an isomorphism
$$H^*(M^*)\simeq H_{n-*}(Y,\partial Y).$$
Under this isomorphism, the connecting homomorphism $\delta$ in the relative
homology exact sequence
$$H_*(Y)\longrightarrow H_*(Y,\partial Y)\stackrel{\delta}{\longrightarrow}
H_{*-1}(\partial Y)$$
is given by
$$\delta(x)={\cal A}(x)\cap\partial Y.$$
\propositionEnd

(Classical references for this are Thom \cite{thom}, Smale \cite{smale},
and Milnor \cite{milnor:h-cobordism}.  A more general statement is proved in
Floer \cite{floer}.  For more novel proofs see Witten \cite{witten:morse},
Helffer-Sj\"{o}strand \cite{helffer-sjostrand}, and Schwarz
\cite{schwarz}.)

\bigskip

\proof
Define $G:C_*(Y)\to M^{n-*}$ by
$$G(\alpha)=\sum_{x\in\mbox{\scriptsize
Crit}}({\cal D}(x)\cdot\alpha)x.$$
Clearly $G$ annihilates $C_*(\partial Y)$ and therefore defines a map
$C_*(Y,\partial Y)\to M^{n-*}$.  Applying $\partial$ to
Lemma~\ref{lem:flowGeometry} and using Lemma~\ref{lem:AscChainMap}, we obtain
$$0=-{\cal A}(G(\partial\alpha))+{\cal A}(dG(\alpha))$$
in $C_*(Y,\partial Y)$.  It follows that $G$ is a chain map.  We
claim that the induced map on homology is the inverse of the map induced by
${\cal A}$.  By definition, $G\circ{\cal A}$ is equal to the identity
on $N_*$.  On the other hand, Lemma~\ref{lem:flowGeometry} asserts that
${\cal F}$ is a chain homotopy between ${\cal A}\circ G$ and the
identity on $C_*(Y,\partial Y)$.

The assertion about the connecting
homomorphism is true more or less by definition.
\qed

We now consider the generalization to circle-valued Morse functions.  Let
the notation be as in the introduction.  Let $C_*(\tilde{X},+\infty)$ be
the complex of locally finite chains in $\tilde{X}$ that are supported away
from the lower end of $\tilde{X}$, i.e.\ for each chain there exists
$R\in{\Bbb R}$ such that the chain is supported in
$\{(x,\lambda)\in\tilde{X} \mid \lambda\ge R\}$.  If $x\in\mbox{Crit}$, let
${\cal A}(t^kx)\subset\tilde{X}$ denote the ascending manifold of a lift
$(x,\lambda)\in\tilde{X}$ of $x$ with $k<\lambda<k+1$.

The following (without ${\cal A}$) was observed by Novikov \cite{novikov}:

\proposition
\label{pro:novikov}
${\cal A}$ is a chain map and induces an isomorphism of $L_{{\Bbb Z}}$
modules
$$H^*(M^*)\simeq H_{n-*}(\tilde{X},+\infty).$$
\propositionEnd

\proof
Same idea as the proof of Proposition~\ref{pro:morseHomology}.
\qed

This proposition has two corollaries which we will also need.

\corollary
\label{cor:isomorphism}
Let $C_*(\tilde{X})$ be the chain complex associated to an equivariant cell
decomposition of $\tilde{X}$, as in the introduction.  Then we have an
isomorphism of $L_{{\Bbb Z}}$-modules
$$H_*\left(C_*(\tilde{X})\bigotimes_{{\Bbb Z}[t,t^{-1}]}
L_{{\Bbb Z}}\right)\simeq H^{n-*}(M^*).$$ 
\corollaryEnd

\proof
There is an obvious chain map
$$C_*(\tilde{X})\bigotimes_{{\Bbb Z}[t,t^{-1}]} L_{{\Bbb Z}}\longrightarrow
C_*(\tilde{X},+\infty).$$
This induces an isomorphism in homology, thanks to the equivalence between
cellular and singular homology.  Now compose this isomorphism with the
isomorphism of Proposition~\ref{pro:novikov}.
\qedNoSkip

\corollary
\label{cor:surjective}
Suppose $H^*(M^*\otimes L_{{\Bbb Q}})=0$.  Then the map
$H_*(\Sigma;{\Bbb Q})\to H_*(\tilde{X};{\Bbb Q})$
induced by the inclusion $p\mapsto (p,0)$ is surjective.
\corollaryEnd

\proof
Let $\alpha\in C_*(\tilde{X};{\Bbb Q})$ be a cycle.  By
Proposition~\ref{pro:novikov} (tensored with $L_{{\Bbb Q}}$)
we can write $\alpha=\partial\beta$ for some $\beta\in
C_*(\tilde{X},+\infty;{\Bbb Q})$.  Observe that
\begin{equation}
\label{eqn:beta}
\partial(\beta\cap\phi^{-1}(-\infty,0]) = \alpha\cap\phi^{-1}(-\infty,0]
+ \beta\cap(\Sigma\times\{0\}).
\end{equation}
Now we can turn this reasoning upside down and apply
Proposition~\ref{pro:novikov} to the Morse function $-\phi$ to find
$\gamma\in C_*(\tilde{X},-\infty;{\Bbb Q})$ with $\alpha=\partial\gamma$.
(The Morse complex for $-\phi$ is acyclic because it is just the dual of
the Morse complex for $\phi$.) We have 
\begin{equation}
\label{eqn:gamma}
\partial(\gamma\cap\phi^{-1}[0,\infty)) = \alpha\cap\phi^{-1}[0,\infty)
- \gamma\cap(\Sigma\times\{0\}).
\end{equation}
Subtracting (\ref{eqn:gamma}) from (\ref{eqn:beta}), we see that
$(\gamma-\beta)\cap(\Sigma\times\{0\})$ is a cycle in $\Sigma\times\{0\}$
homologous to $\alpha$.
\qedNoSkip

\subsection{Torsion}
\label{sec:torsion}

Let $R$ be an integral domain and let
$$C^0\stackrel{d}{\longrightarrow}
C^1\stackrel{d}{\longrightarrow}\cdots\stackrel{d}{\longrightarrow}C^m$$ be
a complex of finitely generated free $R$-modules.  Suppose that $C^*\otimes
Q(R)$ is acyclic.  Using a volume form on $C^i\otimes Q(R)$ coming from a
free basis for $C^i$, we can define the Reideimeister torsion
$$\tau(C\otimes Q(R))\in \frac{Q(R)}{\{\mbox{units of $R$}\}}$$
as in (\ref{eqn:tauDef}).  (Choosing a different free basis multiplies
$\tau$ by a unit.)

Note that if $h:Q(R)\to F$ is an inclusion into a larger field, then
$C\otimes F$ is acyclic if and only if $C\otimes Q(R)$ is, and
$$\tau(C\otimes F)=h(\tau(C\otimes Q(R))).$$
We will denote $\tau(C\otimes Q(R))$ simply by $\tau(C)$.

To compute $\tau(C)$, we will use the following result of Turaev
\cite[\S2.1]{turaev}, generalizing a theorem of Milnor
\cite{milnor:coverings}.  If $E$ is a finitely generated module over a ring
$R$, let $\mbox{Fitt}_1(E)$ be the first Fitting ideal of the module, which
is generated by the determinants of the $n\times n$ minors of the matrix of
relations for a presentation of $E$ with $n$ generators.  This does not
depend on the presentation.  (See e.g.\ Rolfsen \cite{rolfsen}.)  If
greatest common divisors exist in $R$, let $\mbox{ord}(E)$ denote the
greatest common divisor of the elments in $\mbox{Fitt}_1(E)$, which is
well-defined up to units in $R$.

\lemma
\label{lem:torsion}
Suppose $R$ is a Noetherian UFD, and $C$ is as above. Then
$$\tau(C)=\prod_{i=0}^m(\mbox{ord}\;H^i(C))^{(-1)^{i}}$$
up to units of $R$.
\lemmaEnd

In particular, the rings ${\Bbb Z}[t,t^{-1}]$ and $L_{{\Bbb Z}}$ satisfy the
hypothesis on $R$ above.  This lemma is an easy exercise in the special
case when $\mbox{Ker}(d)$ is a free $R$-module with a free complement.

\subsection{Equivalence of Theorems 1.1 and 1.2}
\label{sec:theoremsEquivalent}

This follows from:

\lemma
\begin{description}
\item{(a)}
The leading coefficients of the left and right sides of
Theorem~\ref{thm:refinement} are equal, up to sign.
\item{(b)}
Theorem~\ref{thm:main} is the logarithmic derivative of
Theorem~\ref{thm:refinement}.
\end{description}
\lemmaEnd

\proof
(a)
By definition, the leading coefficient of $\zeta(f)$ is 1, so we need
to check that the leading coefficients of $\tau(M)$ and $\tau(X,\phi)$
agree.  We defined
$$\tau(X,\phi)= \tau\left(C_*(\tilde{X})
\bigotimes_{{\Bbb Z}[t,t^{-1}]}{\Bbb Z}[t,t^{-1}]\right),$$ 
where $C_*(\tilde{X})$ is the complex associated to an equivariant cell
decomposition of $\tilde{X}$.  We can tensor
$C_*(\tilde{X})\bigotimes_{{\Bbb Z}[t,t^{-1}]}{\Bbb Z}[t,t^{-1}]$ with
$L_{{\Bbb Z}}$, and this will 
not affect the torsion (or more precisely will include the torsion into
$L_{{\Bbb Q}}$).  Thus
$$\tau(X,\phi) = \tau\left(C_*(\tilde{X}) \bigotimes_{{\Bbb Z}[t,t^{-1}]}
L_{{\Bbb Z}}\right).$$ 

By Corollary~\ref{cor:isomorphism},
$$H_*\left(C_*(\tilde{X}) \bigotimes_{{\Bbb Z}[t,t^{-1}]}
L_{{\Bbb Z}}\right)\simeq H^{n-*}(M_*).$$ 
From Lemma~\ref{lem:torsion},
we see that the torsion of a complex of $L_{{\Bbb Z}}$-modules depends only
on the homology, up to units in $L_{{\Bbb Z}}$.  Thus
$$\tau(M)=\tau(X,\phi)$$
up to units in $L_{{\Bbb Z}}$.  But a unit in $L_{{\Bbb Z}}$ must have
leading coefficient $\pm 1$, so we are done.

(b)
Observe that
$$\frac{d}{dt}\log\det(1-tA_i)=-\Tr(A_i(1-tA_i)^{-1}).$$
Combining this with (\ref{eqn:zetaDef}), we see that we need
$$\tau(X,\phi)=c\prod_{i=0}^{n}\det(1-tA_i)^{(-1)^{n-i}}$$
for some $c\in{\Bbb R}$, up to units in ${\Bbb Z}[t,t^{-1}]$.  (Note that the
integer $m$ on the right side of Theorem~\ref{thm:main} absorbs this
ambiguity in $\tau(X,\phi)$.)  Now $\tau(X,\phi)$ is the torsion of the
complex $C_*(\tilde{X})\bigotimes_{{\Bbb Z}[t,t^{-1}]}{\Bbb Z}[t,t^{-1}]$,
and the homology of 
this complex is isomorphic to $H_*(\tilde{X})$ as a
${\Bbb Z}[t,t^{-1}]$-module.  So we will show that
$$\mbox{ord}(H_i(\tilde{X}))=c_i\det(1-tA_i)$$
for some $c_i\in{\Bbb Z}$, and then we will be done by Lemma~\ref{lem:torsion}.

We can choose a set $S\subset H_i(\tilde{X})$ which projects to a basis for
$H_i(\tilde{X})/\mbox{Torsion}$ over ${\Bbb Z}$.  ($S$ will be finite thanks to
Corollary~\ref{cor:surjective}.)
We can then choose a finite set $T$ that generates the
torsion part of $H_i(\tilde{X})$ as a ${\Bbb Z}[t,t^{-1}]$-module.  The
matrix of relations for $H_i(\tilde{X})$ is the following:
\vbox{
$$
\;\;\;\;\;\;\;\;\;\;\;S\;\;\;\;\;\;\;\;T\;\;\;T
$$
$$
\begin{array}{c}S \\ T \end{array}
\left(
\begin{array}{ccc}
1-tA_i & 0 & 0\\
? & D & ?
\end{array}
\right)
$$
}
Here the columns represent relations.  The only relations on
$H_i(\tilde{X};{\Bbb Q})$ are $1-tA_i$, and when we lift these to
$H_i(\tilde{X})$ via our choice of $S$, there may be an additional
component in $T$, which is the lower left block of the matrix.  $D$ is a
diagonal matrix of integers asserting that the elements of $T$ are torsion.
The lower right block of the matrix expresses whatever additional relations
the elements of $T$ may satisfy amongst themselves.

Now every minor of this matrix is divisible by $\det(1-tA_i)$, so
$\det(1-tA_i)$ divides $\mbox{ord}(H_i(X))$.  On the other hand one of the
minors is $\det(D)\det(1-tA_i)$, so $\mbox{ord}(H_i(X))$ divides
$\det(D)\det(1-tA_i)$.  Our claim follows.
\qedNoSkip

\sector{Proof of Theorem 1.1}

\subsection{Outline of the proof}
\label{sec:outline}

One of the classical proofs of the Lefschetz fixed point formula on a
manifold $\Sigma$ goes as follows.  We wish to calculate the intersection
number of the graph with the diagonal in $\Sigma\times\Sigma$.  We can
replace the diagonal with a homologous cycle in $C_*(\Sigma)\otimes
C_*(\Sigma)$, and then we are reduced to intersection theory in $\Sigma$.
We will attempt to extend this reasoning to our situation, where the graph
is no longer a cycle.  Since intersection number with a chain with boundary
does not descend to homology, more care is required.

First of all let us assume that $\partial X=\emptyset$;  it is easy to
remove this restriction at the end, in \S\ref{sec:boundary}.

Define a chain $\Gamma$ in $\Sigma\times\Sigma$ with Taylor series
coefficients by
$$\Gamma \eqdef \sum_{k=1}^\infty t^k(\mbox{graph of $f^k$}).$$
Let $\mbox{diag}\subset\Sigma\times\Sigma$ be the diagonal.  Then the first
term in Theorem~\ref{thm:main} is the intersection number
$\Gamma\cdot\mbox{diag}$.  

We wish to replace $\mbox{diag}$ with a homologous cycle in
 $C_*(\Sigma;L_{{\Bbb Q}})^{\otimes 2}$.  However $\Gamma$ is not a
cycle, as its closure has nontrivial boundary involving the ascending and
descending manifolds of critical points.  (We will not distinguish $\Gamma$
from its closure in the notation.)  So we will first find $Z\in
C_*(\Sigma;L_{{\Bbb Q}})^{\otimes 2}$ such that
$$\partial(\Gamma-Z)=0.$$
Such a $Z$ exists by the Eilenberg-Zilber theorem, but is not canonical.
However we will see in \S\ref{sec:closeOffBoundary} that when the rational
Morse complex is acyclic, there is a canonical choice of $Z$, constructed
directly out of the gradient flow.

Next, let $\{e_i\}$ be a set of cycles in $\Sigma$ that represent a basis
for $H_*(\Sigma;{\Bbb Q})$, and let $\{e_i^*\}$ be cycles representing the
(Poincar\'{e}) dual basis, i.e. $e_i\cdot e_j^*=\delta_{ij}$.  Then
$\mbox{diag}-\sum_ie_i\times e_i^*$ is homologous to zero.  It follows that
the intersection number
$$(\Gamma-Z)\cdot(\mbox{diag}-\sum_ie_i\times e_i^*)=0.$$

Direct calculations in \S\ref{sec:calculatingIntersections} will show that
for the natural choice of $Z$ mentioned above,

\lemma
\label{lem:intersections}
\begin{description}
\item{(a)}
$Z\cdot\mbox{diag}  =  (-1)^nt\frac{d}{dt}\log\tau(M)$.
\item{(b)}
$(\Gamma-Z)\cdot\sum_ie_i\times e_i^*  =  \Tr(B:H_*(\Sigma;L_{{\Bbb Q}})\to
H_*(\Sigma;L_{{\Bbb Q}})))$.
\end{description}
\lemmaEnd

Here $B\in\End(C_*(\Sigma;L_{{\Bbb Q}}))$ is a natural chain map constructed in
\S\ref{sec:formalism} out of the gradient flow.  Roughly speaking, $B$ is
$\sum_{k=1}^{\infty}t^kf^k$, plus a correction term that makes it a chain map.

In \S\ref{sec:understandingB} we prove:

\lemma
\label{lem:diagramCommutes}
The diagram
$$
\begin{array}{ccc}
H_*(\Sigma) & \stackrel{B}{\longrightarrow} & H_*(\Sigma) \\
\big\downarrow\imath_* & & \big\downarrow\imath_*\\
H_*(\tilde{X}) & \stackrel{tA(1-tA)^{-1}}{\longrightarrow} & H_*(\tilde{X})
\end{array}
$$
commutes.  (Here all homology is with $L_{{\Bbb Q}}$ coefficients, and
$\imath:\Sigma\to\tilde{X}$ sends $p\mapsto(p,0)$.)
\lemmaEnd

The proof of this lemma is quite natural.  Let $\gamma\in
C_*(\Sigma;{\Bbb Q})$ be a cycle, and suppose the upward gradient flow takes
$\gamma$ around $X$ $k$ times without hitting any critical points.  Then
$\gamma\times\{-k\}$ and $f^k(\gamma)\times\{0\}$ are homologous in
$\tilde{X}$, because their difference is the boundary of the entire
gradient flow between them.  But $\gamma\times\{-k\}$ is the $k^{th}$
downward deck transformation of $\gamma\times\{0\}$, so the $k^{th}$
downward deck transformation of $\gamma\times\{0\}$ is homologous in
$\tilde{X}$ to $f^k(\gamma)\times\{0\}$.  This means that
$\imath_*f^k(\gamma)=A^k(\imath_*\gamma)$.  More generally, if the upward
gradient flow of $\gamma$ hits some critical points, then the gradient flow
no longer gives a homology between $\gamma\times\{-k\}$ and
$f^k(\gamma)\times\{0\}$, because it has additional boundary components
arising from the critical points.  But the extra term in $B$ is exactly
what is needed to cancel these.

By Corollary~\ref{cor:surjective}, $\imath_*$ is surjective, so
$$\Tr(B)=\Tr(tA(1-tA)^{-1})+\Tr(B|\Ker(\imath_*))$$
So if we can show $\Tr(B|\Ker(\imath_*))\in{\Bbb Z}$, we are done.  We can
understand $\Ker(\imath_*)$ in terms of the Morse theory using
Proposition~\ref{pro:morseHomology}.  We then compute the restriction of
$B$ to this kernel, and we use a cheap trick to show that its trace is an
integer: we argue that all the nonconstant terms in $\Tr(B|\Ker(\imath_*))$
vanish a priori, basically because they have the wrong degree.

To carry out the above computations,  we need to develop some formalism.
\S\ref{sec:geometry} proves some simple geometrical facts we need, along
the lines of Lemma~\ref{lem:flowGeometry}, and \S\ref{sec:formalism} encodes
these geometrical facts into algebraic formalism.

\subsection{Geometric observations}
\label{sec:geometry}

Let $Y$ be an oriented manifold and let $\psi:Y\to[0,r]$ be a  generic Morse
function.  Let $Y_0=\psi^{-1}(0)$, $Y_1=\psi^{-1}(r)$, and assume $\partial
Y = Y_0\union Y_1$.  If $x$ is a
critical point, define ${\cal A}(x)$ and ${\cal D}(x)$, and orient
them, as in \S\ref{sec:morseHomology}.  The gradient flow of $\psi$ defines
a diffeomorphism
$$g:Y_0\setminus\bigcup_{x\in\mbox{Crit}}{\cal D}(x)\longrightarrow
Y_1\setminus\bigcup_{x\in\mbox{Crit}}{\cal A}(x).$$

\lemma
\label{lem:geometry}
\begin{description}
\item{(a)}
$\partial(\mbox{graph of $g$}) = \sum_i(-1)^{i}
\sum_{x\in\mbox{\scriptsize Crit}^i} ({\cal D}(x)\cap Y_0) \times
({\cal A}(x)\cap Y_1).$
\item{(b)} If $\alpha$ is a chain in $Y_0$, then
$$
\partial g(\alpha) = g(\partial \alpha)-\sum_{x\in\mbox{\scriptsize
Crit}}({\cal D}(x)\cdot \alpha)({\cal A}(x)\cap Y_1). 
$$
\item{(c)} If $\alpha$ is a chain in $Y_1$, then
$$\partial g^{-1}(\alpha)=g^{-1}(\partial \alpha) +
\sum_{x\in\mbox{\scriptsize Crit}} (\alpha\cdot{\cal A}(x))
({\cal D}(x)\cap Y_0).$$
\end{description}
\lemmaEnd

\proofSketch
(a) is straightforward.  (As in the proof of Lemma~\ref{lem:flowGeometry},
it is easiest to first assume all the critical points are at the same
height and then use induction.)  (b) follows by applying $\partial$ to
Lemma~\ref{lem:flowGeometry}.  (c) is analogous to
(b).
\qedNoSkip

\subsection{Some formalism}
\label{sec:formalism}

If $x\in\mbox{Crit}^i$, let ${\cal A}(t^kx)\subset\tilde{X}$ (resp.\
${\cal D}(t^kx)$) be the ascending (resp.\ descending) manifold of a lift
$(x,\lambda)\in\tilde{X}$ of $x$ with $k < \lambda < k+1$.  Define a map
$$\pi_+:M^*\otimes L_{{\Bbb Z}}\to C_{n-1-*}(\Sigma;
L_{{\Bbb Z}})$$
by setting
$$\pi_+(x) \eqdef \sum_{j=-\infty}^\infty
t^{j-1}{\cal A}(x)\cap(\Sigma\times\{j\})$$
for all $x\in M^*$.
(Here we are identifying $\Sigma\times\{j\}$ with $\Sigma$.)

In English,
$\pi_+$ of a critical point is the set of all points in $\Sigma$ such that
downward gradient flow converges to the critical point, multiplied by
$t^k$, where $k$ is the number of times the gradient flow crosses through
$\Sigma$ before reaching the critical point.

We define
$$\pi_-:M^*\otimes L_{{\Bbb Z}}\to C_{*-1}(\Sigma;L_{{\Bbb Z}})$$
similarly; for $x\in M^*$, let
$$\pi_-(x) \eqdef \sum_jt^j{\cal D}(x)\cap(\Sigma\times\{-j\}).$$

\lemma
\label{lem:boundaryGraph}
$\partial\Gamma = t\sum_{i}{(-1)^i} \sum_{y\in\mbox{\scriptsize Crit}^i}
\pi_-(y)\times\pi_+(y).$ 
\lemmaEnd

\proof
We use Lemma~\ref{lem:geometry}(a).  To
calculate the $t^k$ term, let
$$Y=\{(p,\lambda)\in \tilde{X}\mid 0\le \lambda\le
k\}\subset\tilde{X}.$$
So in the notation of \S\ref{sec:geometry}, $Y_0=Y_1=\Sigma$ and $g=f^k$.

If $y\in\mbox{Crit}$, let $\phi(y)_0\in(0,1)$ be a representative of
$\phi(y)\in{\Bbb R}/{\Bbb Z}$.  Then the critical points of $g$ are of the
form $(y,\phi(y)_0+j)$ for $y\in\mbox{Crit}$ and $j=0,\ldots,k-1$.  The
intersection of the
descending manifold of such a critical point with $\phi^{-1}(0)$ is the
$t^j$ term of $\pi_-(y)$.  The intersection of the ascending manifold with
$\phi^{-1}(k)$ is the $t^{k-j-1}$ term of $\pi_+(y)$.  So by
Lemma~\ref{lem:geometry}(a),
$$
\partial(\mbox{graph of $f^k$})  =  \sum_{i}{(-1)^i}
\sum_{y\in\mbox{\scriptsize Crit}^i} \sum_{j=0}^{k-1} (\mbox{$t^j$ term of
$\pi_-(y)$}) \times (\mbox{$t^{k-j-1}$ term of $\pi_+(y)$}).
$$
This proves the $t^k$ term of the lemma.
\qed

In the sequel, we will omit the details when applying
Lemma~\ref{lem:geometry} in a straightforward way as above.

\lemma
\label{lem:piChainMap}
If $d$ is the differential in the Morse complex and $d^*$ is its adjoint
with respect to the natural inner product $\langle\;,\;\rangle$, then
\begin{description}
\item{(a)} $\partial\pi_+=-\pi_+d$,
\item{(b)} $\partial\pi_-=(-1)^{i-1}\pi_-d^*$ on $M^i$.
\end{description}
\lemmaEnd

\proof
Part (a) follows from Lemma~\ref{lem:geometry}(b).  Part (b) follows
similarly from Lemma~\ref{lem:geometry}(c).  The sign here arises when we
switch from intersections in $X$ (as in Lemma~\ref{lem:geometry}) to
intersections in $\Sigma$ (in the definition of $d$ in
\S\ref{sec:morseHomology}).
\qed

The following will be used in \S\ref{sec:calculatingIntersections}.
Let $f$ be the partially defined endomorphism of $\Sigma$ from the
introduction.  Define an endomorphism $f^+$ of $C_*(\Sigma)\otimes
L_{{\Bbb Z}}$ by
$$
f^+  \eqdef  (1-tf)^{-1}.
$$

It turns out that $f^+$ is not a chain map.  To understand this, define a
map $\xi:C_*(\Sigma)\to M^{n-*}$ by requiring that
$$\langle x,\xi(\alpha)\rangle \eqdef \pi_-(x)\cdot\alpha$$
for all $\alpha\in C_*(\Sigma)$ and $x\in M^{n-*}$.  (Note that $\xi$ is
an analogue of the map $G$ in the proof of
Proposition~\ref{pro:morseHomology}.)

\lemma
\label{lem:fPlusNotChainMap}
If $\alpha\subset\Sigma$ is a generic chain, then
$$\partial f^+(\alpha) = f^+(\partial\alpha) - t\pi_+\xi(\alpha). $$
\lemmaEnd

\proof
This follows from Lemma~\ref{lem:geometry}(b).
\qed

We can add a correction term to $f^+$, using $\xi$, to make it a chain map.
(This will arise naturally in \S\ref{sec:calculatingIntersections}.)
Namely, define an endomorphism
$B$ of $C_*(\Sigma;L_{{\Bbb Q}})$ by
$$B \eqdef f^+-1-t\pi_+d^*\Delta^{-1}\xi.$$
($\Delta$ is invertible thanks to our assumption that the rational Morse
complex is acyclic.)

\lemma
\label{lem:xiBChainMap}
\begin{description}
\item{(a)} $d\xi=\xi\partial$.
\item{(b)} $\partial B=B\partial$.
\end{description}
\lemmaEnd

\proof
(a)
Let $\alpha\in C_i(\Sigma)$ and $x\in M^{n-i+1}$. Then
\begin{eqnarray*}
\langle d\xi(\alpha),x\rangle &=&\pi_-d^*x\cdot\alpha\\
&=&(-1)^{n-i}\partial\pi_-x\cdot\alpha\\
&=&\pi_-x\cdot\partial\alpha\\
&=&\langle\xi(\partial\alpha),x\rangle.
\end{eqnarray*}
(The sign in the third line is tricky, because $\partial$ is not quite a
signed derivation with respect to intersections.  Also, an alternate proof
of (a) can be given by mimicking the proof that $G$ is a chain map in
\S\ref{sec:morseHomology}.)

(b)
By Lemma~\ref{lem:fPlusNotChainMap},
\begin{equation}
\label{eqn:xiB1}
\partial B =
f^+\partial-t\pi_+\xi-\partial-t\partial\pi_+d^*\Delta^{-1}\xi.
\end{equation}
We can write
\begin{equation}
\label{eqn:xiB2}
t\pi_+\xi=t\pi_+dd^*\Delta^{-1}\xi+t\pi_+d^*\Delta^{-1}d\xi.
\end{equation}
By Lemma~\ref{lem:piChainMap}(a),
$$t\pi_+dd^*\Delta^{-1}\xi=-t\partial\pi_+d^*\Delta^{-1}\xi.$$
By part (a),
$$t\pi_+d^*\Delta^{-1}d\xi=t\pi_+d^*\Delta^{-1}\xi\partial.$$
Substitute the above two equations into (\ref{eqn:xiB2}), and
put the result into (\ref{eqn:xiB1}). 
\qedNoSkip

\subsection{Closing off the boundary of the graph}
\label{sec:closeOffBoundary}

As explained in the outline, we now want to find $Z\in
C_*(\Sigma;L_{{\Bbb Q}})^{\otimes 2}$ with $\partial
Z=\partial\Gamma$.  Let $P$ be the composition
$$\Hom(M^i,M^j)\stackrel{\rho}{\longrightarrow} M^j\otimes
M^i\stackrel{\pi_-\otimes\pi_+}{\longrightarrow}
C_*(\Sigma;L_{{\Bbb Q}})^{\otimes 2}.$$
Here $\rho$ is the canonical isomorphism given by the inner product
$\langle\;,\;\rangle$.
Our ansatz will be $Z=P(W)$ for some $W\in\Hom(M^*,M^{*+1})$.

\lemma
\label{lem:boundaryP}
Let $W=\sum_{i=0}^{n-1}W^i$ with $W^i\in\Hom(M^i,M^{i+1})$. Then
$$\partial P(W)=\sum_{i=0}^{n-1}(-1)^iP(d^*W_i+W_{i-1}d^*).$$
(Here we interpret $W_{-1}=0$.)
\lemmaEnd

\proof
By Lemma~\ref{lem:piChainMap},
$$\partial P(W)=\sum_i((-1)^i\pi_-d^*\otimes
\pi_+-(-1)^i\pi_-\otimes\pi_+d)\rho(W_i).$$
(The factor of $(-1)^i$ on the right arises because $\partial(a\times
b)=\partial a\times b + (-1)^{\dim(a)}a\times\partial b$, and
$\pi_-M^{i+1}$ has dimension $i$.)  Now use the facts $(d^*\otimes
1)\rho(W_i)=\rho(d^*W_i)$ and $(1\otimes d)\rho(W_i)=\rho(W_id^*)$.
\qed

In this notation, Lemma~\ref{lem:boundaryGraph} says that
$$\partial\Gamma=t\sum_{i=0}^{n-1}(-1)^iP(1|M^i).$$
So by Lemma~\ref{lem:boundaryP}, $\partial Z=\partial\Gamma$ if and only if
$$d^*W+Wd^*=t.$$
Such a $W$ exists if and only if the rational Morse complex is acyclic,
which we assumed to be true.   The natural choice, which we will adopt, is
$$W\eqdef t\Delta^{-1}d.$$

\subsection{Calculating intersection numbers}
\label{sec:calculatingIntersections}

We will now prove Lemma~\ref{lem:intersections}.
In these calculations, it is convenient to choose a basis
$\{x^i_j\}$ for $\Ker(d^*|M^i)$ with $\|x^i_j\|=1$ and $\Delta
x^i_j=\lambda^i_jx^i_j$.  (To find such a basis we may have to extend to
coefficients in the algebraically closed field $L_{{\Bbb C}}$, which causes
no problems.)

\lemma
\label{lem:torsionEigenvalues}
$\tau(M) = \prod_{i=0}^{n-1} \left( \prod_j\sqrt{\lambda^i_j}
\right)^{(-1)^{i-1}}.$
\lemmaEnd

\proof
Take $\omega_i=\bigwedge_jx^i_j$ in (\ref{eqn:tauDef}).
\qedNoSkip

\lemma
\label{lem:piMinusIntersectPiPlus}
Let $x\in M^*$ and $y\in M^{*+1}$.  Then
$$\pi_-y\cdot\pi_+x=\frac{d}{dt}\langle dx,y\rangle - \left\langle
d\left(\frac{d}{dt}x\right),y\right\rangle - \left\langle dx,
\frac{d}{dt}y\right\rangle.$$
\lemmaEnd

\proof
If $x,y\in\mbox{Crit}$, then the two rightmost terms are zero, and the
formula follows directly from the definitions.  The general case follows by
expanding $x$ and $y$ in powers of $t$.
\qed

\proofOf{Lemma~\ref{lem:intersections}}
(a) We have
$$Z=t\sum_{i,j}(\lambda^i_j)^{-1}\pi_-(dx^i_j)\times\pi_+(x^i_j).$$
If $\alpha,\beta$ are two chains of complementary dimension then
$(\alpha\times\beta)\cdot\mbox{diag}=(-1)^{\dim(\beta)}\alpha\cdot\beta$.
Thus 
\begin{equation}
\label{eqn:ZCdotDiag}
Z\cdot\mbox{diag} =
t\sum_{i,j}(-1)^{n-1-i}(\lambda^i_j)^{-1}\pi_-(dx^i_j)\cdot\pi_+(x^i_j).
\end{equation}
Writing $x=x^i_j$, Lemma~\ref{lem:piMinusIntersectPiPlus} gives
$$\pi_-(dx)\cdot\pi_+(x)=\frac{d}{dt}\|dx\|^2 - \left\langle
d\left(\frac{d}{dt}x\right),dx\right\rangle - \left\langle dx,
\frac{d}{dt}dx\right\rangle.$$
Thanks to the special properties of $x=x^i_j$, the middle term on the right
vanishes:
\begin{eqnarray*}
\left\langle  d\left(\frac{d}{dt}x^i_j\right),dx^i_j\right\rangle & = & 
\left\langle  \frac{d}{dt}x^i_j,d^*dx^i_j\right\rangle\\
&=&\lambda^i_j\left\langle \frac{d}{dt}x^i_j,x^i_j\right\rangle\\
&=&\frac{\lambda^i_j}{2}\frac{d}{dt}\|x^i_j\|^2\\
&=&0.
\end{eqnarray*}
Thus
$$\pi_-(dx^i_j)\cdot\pi_+(x^i_j)  =  \frac{1}{2}\frac{d}{dt}\|dx^i_j\|^2
 = \frac{1}{2}\frac{d}{dt}\lambda^i_j.$$
Substituting this into (\ref{eqn:ZCdotDiag}) and comparing with
Lemma~\ref{lem:torsionEigenvalues} proves (a).

(b)
First of all,
\begin{equation}
\label{eqn:GammaCdotEk}
\Gamma\cdot\sum_ke_k\times e_k^* = \sum_{n=1}^\infty
t^n\sum_i(-1)^{\dim(e_k)}f^n(e_k)\cdot e_k^*.
\end{equation}
Second,
\begin{eqnarray}
Z\cdot\sum_ke_k\times e_k^* & = &
t\sum_{i,j,k}(-1)^{\dim(e_k)}((\lambda^i_j)^{-1}\pi_-dx^i_j\cdot
e_k)(\pi_+(x^i_j)\cdot e_k^*)\nonumber\\ 
&=&t\sum_{i,j,k}(-1)^{\dim(e_k)}(\lambda^i_j)^{-1}\langle
\xi(e_k),dx^i_j\rangle(\pi_+(x^i_j)\cdot e_k^*). \label{eqn:ZCdotEk}
\end{eqnarray}
By Lemma~\ref{lem:xiBChainMap}(a), $d\xi(e_k)=0$, so we can write
$$\xi(e_k)=dd^*\Delta^{-1}\xi(e_k).$$
Then
\begin{eqnarray*}
\langle \xi(e_k),dx^i_j\rangle & = & \langle
d^*\Delta^{-1}\xi(e_k),d^*dx^i_j\rangle\\ 
&=&\lambda^i_j\langle d^*\Delta^{-1}\xi(e_k),x^i_j\rangle.
\end{eqnarray*}
Putting this into (\ref{eqn:ZCdotEk}) gives
$$
Z\cdot\sum_ke_k\times e_k^* =
t\sum_k(-1)^{\dim(e_k)}(\pi_+d^*\Delta^{-1}\xi(e_k))\cdot e_k^*.$$ 
Subtracting this from (\ref{eqn:GammaCdotEk}) gives
$$Z\cdot(\Gamma-\sum_ke_k\times e_k^*)=\sum_k(-1)^{\dim(e_k)}B(e_k)\cdot e_k^*.
\qedInEquation

\subsection{Understanding $B$}
\label{sec:understandingB}

We will now prove Lemma~\ref{lem:diagramCommutes}.
To do this we need a bit more formalism. If $\alpha$ is a
chain in $\Sigma$, define
$${\cal F}(t^k\alpha) \eqdef
{\cal F}(\alpha\times\{k\})\subset\tilde{X}.$$
If $\alpha\in C_*(\Sigma;L)$ and $k\in{\Bbb Z}$, define
$${\cal F}_k(\alpha) \eqdef \{(x,\lambda)\in{\cal F}(\alpha) \mid
\lambda\le k\}.$$
If $x\in M^*$, define
$${\cal A}_k(x) \eqdef \{(x,\lambda)\in{\cal A}(x) \mid \lambda\le k\}.$$

\lemma
\label{lem:moreFlowGeometry}
\begin{description}
\item{(a)} If $\gamma\in C_*(\Sigma;{\Bbb Q})$ is a cycle and $k>0$, then
$$\partial{\cal F}_k(\gamma) = f^k(\gamma)\times\{k\}-\gamma\times\{0\}
+ {\cal A}_k(\xi(\gamma)).$$
\item{(b)} If $x\in M^*$ then
$$\partial{\cal A}_k(x)={\cal A}_k(dx)+(\pi_+x)^{k-1}\times{k},$$
where $(\pi_+x)^{k-1}$ is the coefficient of $t^{k-1}$ in $\pi_+x$.
\end{description}
\lemmaEnd

\proof
(a) follows directly from Lemma~\ref{lem:flowGeometry}.  (b) follows from
Lemma~\ref{lem:AscChainMap}.
\qed

\proofOf{Lemma~\ref{lem:diagramCommutes}}
Let $\gamma\in C_*(\Sigma;{\Bbb Q})$ be a cycle.  We need to show:
\begin{description}
\item{(a)} If $k>0$, then $(B\gamma)^k\times\{k\}$ is homologous to
$\gamma\times\{0\}$ in $\tilde{X}$ (where $(B\gamma)^k$ is the $t^k$
coefficient in $B(\gamma)$).
\item{(b)} If $k\le 0$, then $(B\gamma)^k\times\{k\}$ is nullhomologous in
$\tilde{X}$.
\end{description}

Suppose $k>0$.  By Lemma~\ref{lem:xiBChainMap}(a),
$d\xi(\gamma)=0$, so we can write
$$\xi(\gamma)=d(d^*\Delta^{-1}\xi \gamma).$$
Putting $x=d^*\Delta^{-1}\xi \gamma$ into
Lemma~\ref{lem:moreFlowGeometry}(b) gives
$$\partial{\cal A}_k(d^*\Delta^{-1}\xi \gamma) = {\cal A}_k(\xi
\gamma)+(\pi_+d^*\Delta^{-1}\xi \gamma)^{k-1}\times\{k\}.$$ 
Subtracting this from Lemma~\ref{lem:moreFlowGeometry}(a) gives
\begin{eqnarray*}
\partial(\mbox{something}) & = &
f^k(\gamma)\times\{k\}-(\pi_+d^*\Delta^{-1}\xi
\gamma)^{k-1}\times\{k\}-\gamma\times\{0\}\\ 
&=&(B\gamma)^k\times\{k\}-\gamma\times\{0\}.
\end{eqnarray*}
This proves (a).

If $k\le 0$, then ${\cal F}_k(\xi(\gamma))=0$, since
$$\xi(\gamma)=\sum_{y\in\mbox{\scriptsize Crit}}(\pi_-(y)\cdot\gamma)y$$
is a Taylor series.  Then Lemma~\ref{lem:moreFlowGeometry}(b) gives
$$\partial{\cal A}_k(d^*\Delta^{-1}\xi\gamma) =
(\pi_+d^*\Delta^{-1}\xi\gamma)^{k-1}\times\{k\}.$$
This implies (b).
\qed

Let $V_+\subset H_*(\Sigma;{\Bbb Q})$ be the subspace generated by cycles of
the form $(\pi_+x)^0$, where $x\in M^*$ and $(dx)^{\le 0}=0$.  (Here
$(\pi_+x)^0$ denotes the constant coefficient of $\pi_+x$, and $(dx)^{\le
0}$ the portion of $dx$ containing nonpositive powers of $t$.)  Similarly,
let $V_-\subset H_*(\Sigma;{\Bbb Q})$ be the subspace generated by cycles of
the form $(\pi_-y)^0$, where $y\in M^*$ and $(d^*y)^{\le 0}=0$.

\lemma
\label{lem:kernel}
\begin{description}
\item{(a)}
$\Ker(\imath_*:H_*(\Sigma;{\Bbb Q}) \to H_*(\tilde{X};{\Bbb Q})) = V_+
\oplus V_-.$
\item{(b)}
Our assumption that $H^*(M^*\otimes L_{{\Bbb Q}})=0$ implies $V_+\intt
V_-=\{0\}$. 
\end{description}
\lemmaEnd

\proof
Define
\begin{eqnarray*}
\tilde{X}^+ & \eqdef & \{(x,\lambda)\in\tilde{X}\mid\lambda\ge 0\},\\
\tilde{X}^- & \eqdef & \{(x,\lambda)\in\tilde{X}\mid\lambda\le 0\}.
\end{eqnarray*}
The relative homology exact sequence
$$H_{k+1}(\tilde{X}^-,\Sigma) \stackrel{\delta}{\longrightarrow}
H_k(\Sigma) \longrightarrow H_k(\tilde{X}^-).$$
and Proposition~\ref{pro:morseHomology} imply that the kernel of the map
$H_*(\Sigma) \to H_*(\tilde{X}^-)$ is $V_+$. (Here we are identifying
$\Sigma$ with $\Sigma\times\{0\}\subset\tilde{X}$, and all homology is with
rational coefficients.) Also
$H_*(\tilde{X},\tilde{X}^-)\simeq H_*(\tilde{X}^+,\Sigma)$ by excision,
and the connecting homomorphism $\delta$ in the exact sequence
$$H_{k+1}(\tilde{X},\tilde{X}^-) \stackrel{\delta}{\longrightarrow}
H_k(\tilde{X}^-) \longrightarrow H_k(\tilde{X}).$$
sends this to $V_-$.  This proves (a).

To prove (b), suppose $u\in V_+\intt V_-$.  Write
$u=(\pi_+x)^0=(\pi_-y)^0$.  Let $v\in H_{k+1}(\tilde{X})$ be the cycle
obtained by gluing together the upward gradient flow of $x$ (up to
$\Sigma$) and the downward gradient flow of $y$ (down to $\Sigma$).  Note
that $u$ is the image of $v$ under the connecting homomorphism $\delta$ in
the Mayer-Vietoris sequence
$$H_{k+1}(\tilde{X}^-)\oplus H_{k+1}(\tilde{X}^+)\longrightarrow
H_{k+1}(\tilde{X})\stackrel{\delta}{\longrightarrow} H_k(\Sigma).$$
By Corollary~\ref{cor:surjective}, $v$ is in the image of
$\imath_*:H_{k+1}(\Sigma)\to H_{k+1}(\tilde{X})$.  But $\delta \imath_*=0$,
so $u=0$.
\qed

We will now compute $B|\Ker(\imath:H_*(\Sigma;{\Bbb Q})\to
H_*(\tilde{X};{\Bbb Q}))$.  Let $R$ be an operator that sends $t^n$ to
$t^{-n}$.

\lemma
\label{lem:traceBOnKernel}
\begin{description}
\item{(a)} If $(dx)^{\le 0}=0$ then
$B((\pi_+x)^0)=-(\pi_+x)^{\le 0}.$
\item{(b)} If $(d^*y)^{\le 0}=0$ then
$B((\pi_-y)^0)=R((\pi_-y)^{<0}).$
\end{description}
\lemmaEnd

\proof
We might as well assume that $x^{>0}=0$.  Then
$$f^+((\pi_+x)^0)=(\pi_+x)^{\ge 0}.$$
From the definitions of $\xi$ and $d$, we have
$$
\xi((\pi_+x)^0) =  t^{-1}(dx)^{>0} = t^{-1}dx.
$$
Then
\begin{eqnarray*}
t\pi_+d^*\Delta^{-1}\xi((\pi_+x)^0) & = & \pi_+x-\pi_+dd^*\Delta^{-1}x\\
&=&\pi_+x + (\mbox{nullhomologous cycle})
\end{eqnarray*}
by Lemma~\ref{lem:piChainMap}.  Putting this into the definition of $B$
proves (a).

To prove (b), let $\gamma$ be a perturbation of $(\pi_-y)^0$.  (We need to
do this because $(\pi_-y)^0$ does not intersect the descending manifolds of
the critical points in $y$ transversely, so $(\pi_-y)^0$ is not a generic
chain on which $f^+$ is defined.) We create $\gamma$ by replacing
$(\pi_y)^0$ with the intersection of $\Sigma\times\{0\}$ and the descending
flow of a sum of small spheres linking the ascending manifolds of the critical
points in $y$ at generic points.  For any positive integer $k$, we can
choose these spheres to be small enough that
$$f^+(\gamma)  =  R((\pi_-y)^{\le0})+(\mbox{nullhomologous cycle})+O(t^k).$$
Similarly
\begin{eqnarray*}
\xi(\gamma) & = &  \pm R((d^*y)^{\le 0})+O(t^k)\\
 & = & O(t^k).
\end{eqnarray*}
We can then complete the proof as in part (a).
\qedNoSkip

\lemma
$\Tr(B|\Ker(\imath_*))\in{\Bbb Z}$.
\lemmaEnd

\proof
By Lemma~\ref{lem:kernel},
$$\Tr(B|\Ker(\imath_*))=\Tr(B|V_+)+\Tr(B|V_-).$$
Now here is the cheap trick.  Since
$$(\Gamma-Z)\cdot(\mbox{diag}-\sum_ie_i\times e_i^*)=0,$$
the combination of Lemmas~\ref{lem:intersections},
Lemma~\ref{lem:diagramCommutes}, and
Corollary~\ref{cor:surjective} implies that
$$\Tr(B|\Ker(\imath_*)) = \Gamma\cdot\mbox{diag}
-(-1)^nt\frac{d}{dt}\log\tau(M)-\Tr(tA(1-tA)^{-1}).$$ 
In particular, we see that $\Tr(B|\Ker(\imath_*))$ is a Taylor series.  For
$\Gamma\cdot\mbox{diag}$ and $\Tr(tA(1-tA)^{-1})$ contain only positive
powers of $t$ by definition, and it is easy to see that
$t\frac{d}{dt}\log$ of a Laurent series is a Taylor series.

It follows from Lemma~\ref{lem:traceBOnKernel} that $\Tr(B|V_+)\in{\Bbb Z}$,
since all the negative degree terms must vanish. We also see from
Lemma~\ref{lem:traceBOnKernel} that the coefficients of $\Tr(B|V_-)$ are
exactly minus what the nonconstant coefficients of $\Tr(B|V_+)$ would be if
we inverted the Morse function (with appropriate new orientations on the
ascending and descending manifolds).  So $\Tr(B|V_-)=0$.
\qed

This completes the proof of Theorems~\ref{thm:main} and
\ref{thm:refinement} in the case $\partial X=\emptyset$, as explained in
\S\ref{sec:outline}.

\subsection{Extension to manifolds with boundary}
\label{sec:boundary}

We wil now deduce Theorem~\ref{thm:refinement} for manifolds with boundary.
Let $L(X,\phi)$ denote the left hand side of Theorem~\ref{thm:refinement}.

\lemma
\label{lem:switchOrientation}
Let $-X$ be $X$ with the opposite orientation. Then
\begin{eqnarray*}
L(X,\phi) & = & L(-X,\phi),\\
\tau(X,\phi) & = & \tau(-X,\phi),
\end{eqnarray*}
modulo the sign ambiguity in $L$ and the $\pm t^k$ ambiguity in $\tau$.
\lemmaEnd

\proof
Suppose we change the orientation of $X$.  The zeta function is not
affected.  We can switch the orientation of the descending manifold of each
critical point to restore the condition that the descending and ascending
manifolds of a critical point have intersection number $+1$ at the critical
point.  After we do this, the differential in the Morse complex is exactly
the same as it was before, so $\tau(M)$ is not affected.  The complex out
of which $\tau(X,\phi)$ is defined is not changed either.
\qedNoSkip

\lemma
\label{lem:gluing}
Let $X$ be closed, and let $Y\subset X$ be a hypersurface separating $X$
into $X_1$ and $X_2$.  Let $\phi:X\to S^1$ be a Morse function such that
$\mbox{grad}(\phi)$ is parallel to $Y$ and nonzero on $Y$.  Assume that the
rational Morse complexes for $X_1$ and $X_2$ are acyclic.  Then the
rational Morse complex for $X$ is acyclic, and
\begin{description}
\item{(a)}
$L(X,\phi)=\frac{L(X_1,\phi)L(X_2,\phi)}{L(Y,\phi)},$
\item{(b)}
$\tau(X,\phi)=\frac{\tau(X_1,\phi)\tau(X_2,\phi)}{\tau(Y,\phi)}.$
\end{description}
\lemmaEnd

\proof
We have
$$M^*(X)=M^*(X_1)\oplus M^*(X_2),$$
where $M^*(X)$ is the Morse complex for $X$, since our assumption on $\phi$
assures that gradient flow lines can never cross $Y$, and there are no
critical points on $Y$. We then trivially obtain
$$\tau(M(X))=\tau(M(X_1))\tau(M(X_2)).$$
Since $\tau(M(Y))=1$, we can rewrite this as
$$\tau(M(X))=\frac{\tau(M(X_1))\tau(M(X_2))}{\tau(M(Y))}.$$
Also every closed orbit lies in either $X_1$ or $X_2$, so we have
$$\zeta(X)=\frac{\zeta(X_1)\zeta(X_2)}{\zeta(Y)}.$$
(We have to divide by $\zeta(Y)$ to avoid counting the closed orbits in $Y$
twice.)  This proves (a).

Next, observe that
$$\tilde{X}=\tilde{X_1}\bigcup_{\tilde{Y}}\tilde{X_2}.$$
Thus we have a short exact sequence of complexes of
${\Bbb Z}[t,t^{-1}]$-modules
$$0\to C_*(\tilde{Y})\to C_*(\tilde{X_1})\oplus C_*(\tilde{X_2})\to
C_*(\tilde{X})\to 0.$$
Now (b) follows from the product formula for torsion \cite[Thm.\
3.1]{milnor:torsion}.
\qed

\proofOf{Theorem~\ref{thm:refinement} when $\partial X\neq\emptyset$}
Given an oriented manifold $X$ with boundary, we can form the double
$2X=X\union_{\partial X}(-X)$.  Since $\mbox{grad}(\phi)$ is parallel to
$\partial X$, we can extend $\phi$ in the obvious way to $2X$, as a $C^1$
function.  Using Lemmas~\ref{lem:switchOrientation} and \ref{lem:gluing},
and applying Theorem~\ref{thm:refinement} to the closed manifolds $\partial
X$ and $2X$, we have
\begin{eqnarray*}
L(X,\phi)^2 & = & L(\partial X,\phi)L(2X,\phi) \\
& = & \tau(\partial X,\phi)\tau(2X,\phi)\\
&=&\tau(X,\phi)^2.
\end{eqnarray*}
Since Theorem~\ref{thm:refinement} is not sensitive to signs, we are done.
\qedNoSkip

\sector{Seiberg-Witten invariants of 3-manifolds}

From now on we assume that $X$ is a closed oriented 3-manifold with
$b^1>0$.  In \S\ref{sec:conjecture}, we define a possible analogue of the
Gromov invariant for $X$, and we conjecture that this is equal to the
Seiberg-Witten invariant.  In \S\ref{sec:meng-taubes} we use
Theorem~\ref{thm:refinement} to show that this conjecture implies the
Meng-Taubes formula (for closed manifolds, up to signs) relating the
Seiberg-Witten invariant to Milnor torsion.

\subsection{An analogue of the Gromov invariant}
\label{sec:conjecture}

Let $\phi:X\to S^1$ be a generic Morse function, and let
$\Sigma=\phi^{-1}(0)$ as usual.  Let $\eta=d\phi$.  Assume that $\phi$ has
no index 0 or 3 critical points.  In particular this implies that the
homology class of $\eta$ is nontrivial.

Let $X'=X\setminus\mbox{Crit}$.  Let $H(\eta)$ denote the set of $\alpha\in
H_1(X',\partial X')$ whose boundary is the sum of the index 2 critical
points minus the sum of the index 1 critical points.  The Gromov invariant
counts unions of flow lines and closed orbits of the gradient flow of
$\phi$ whose homology class is in $H(\eta)$.  We then need some way of
identifying $H(\eta)$ with the set ${\cal S}$ of $Spin^c$-structures.

Let us take care of this last point first.  Given a $Spin^c$ structure $S$
on $X$, let $E\subset S|X'$ be the $-i$ eigenspace of Clifford
multiplication by $\eta/|\eta|$.

\lemma
The map that sends a $Spin^c$ structure $S$ to the Poincar\'{e}-Lefschetz
dual of $c_1(E)$ in $H_1(X',\partial X')$ defines an isomorphism of
$H^2(X;{\Bbb Z})$-torsors
$$\Psi_\eta:{\cal S}\to H(\eta).$$
\lemmaEnd

\proof
Given $S$, the Poincar\'{e} dual of $c_1(E)$ lives in $H(\eta)$ because
$S$ can be trivialized in
a neighborhood of a critical point, but on a sphere around the critical
point, $\eta/|\eta|$ defines a map $S^2\to S^2$ of degree $\pm 1$, so
$E|S^2$ is the Hopf line bundle or its inverse.

Given $E$, we can recover $S$ as follows. Let $K^{-1}$ denote the kernel of
$\eta:TX'\to{\Bbb R}$.  This inherits a complex structure from the metric
and orientation on $X$.  Define
$$S\eqdef E\oplus (K^{-1}\otimes E).$$
The Clifford action is as follows.  If $v\in TX$ is dual to $\eta$, then
$$c(v)\eqdef |v|
\left(\begin{array}{cc}-i & o \\ 0 & i\end{array}\right).$$
If $v\in TX$ is annihilated by $\eta$, i.e.\ $v\in K^{-1}$, then for $e\in E$,
\begin{eqnarray*}
c(v)e & \eqdef & v\otimes e,\\
c(v)(v\otimes e) & \eqdef & -|v|^2e.
\end{eqnarray*}

Every $Spin^c$ structure $S$ must be of this form, because if $E_{\pm}$ is
the $\pm i$ eigenspace of Clifford multiplication by $\eta/|\eta|$, then
Clifford multiplication by $K^{-1}$ defines an isomorphism
$$K^{-1}=\Hom(E_-,E_+),$$
so $E_+=K^{-1}\otimes E_-$, etc.

The $H^2(X;{\Bbb Z})$-torsor structure on $H(\eta)$ is as follows: any two
elements of $H(\eta)$ differ by an element of $H_1(X',\partial X')$ which
is annihilated by $\delta:H_1(X',\partial X')\to H_0(\partial X')$, and
hence extends to an element of $H_1(X) = H^2(X;{\Bbb Z})$.  This clearly
corresponds to the $H^2(X;{\Bbb Z})$ action on ${\cal S}$.
\qed

Now choose orientations on the ascending and descending manifolds of the
critical points, and choose orientations on the vector spaces
${\Bbb Q}^{\mbox{\scriptsize Crit}^1}$ and ${\Bbb Q}^{\mbox{\scriptsize
Crit}^2}$.  Note that $\mbox{Crit}^1$ and $\mbox{Crit}^2$ have the same
cardinality, because $\chi(X)=0$.

Let $\Lambda$ be the Novikov ring of $H_1(X',\partial X')$ with respect to
intersection with $\Sigma$.  This is the ring of functions
$\lambda:H_1(X',\partial X')\to{\Bbb Z}$ such that for any $k\in{\Bbb Z}$,
the set
$$\{\alpha\in H_1(X',\partial X')\mid\lambda(\alpha)\neq 0,\;
\alpha\cdot\Sigma<k\}$$
is finite.  The multiplication is given by the convolution product
$$(\lambda_1\lambda_2)(\alpha) \eqdef
\sum_{\beta}\lambda_1(\beta)\lambda_2(\alpha-\beta).$$
This is a generalization of the ring of Laurent series.  (See
\cite{hofer-salamon} for more discussion.)  We will write elements of
$\Lambda$ like Laurent series, in the form $\sum_i f(\alpha_i)\alpha_i$.

Define a map
$$P:{\Bbb Q}^{\mbox{\scriptsize Crit}^1} \to
{\Bbb Q}^{\mbox{\scriptsize Crit}^2} \otimes \Lambda$$
as follows.  If $x\in\mbox{Crit}^1$ and $y\in\mbox{Crit}^2$, let
${\cal P}(x,y)$ be the set of flow lines from $x$ to $y$ (of the flow dual
to $\eta$), with the orientation induced by $\eta$.  Given
$\gamma\in{\cal P}(x,y)$, let $\epsilon(\gamma)$ be the sign of the
intersection of the descending manifold of $y$ and the ascending manifold
of $x$ (in a local slice orthogonal to $\eta$).  For $x\in\mbox{Crit}^1$,
define
$$P(x) \eqdef \sum_{y\in\mbox{\scriptsize Crit}^2} y \otimes
\sum_{\gamma\in{\cal P}(x,y)} \epsilon(\gamma) [\gamma].$$
Here $[\gamma]$ denotes the homology class of $\gamma$.

Let ${\cal O}$ be the set of closed orbits of the flow dual to $\eta$.  For
$\gamma\in{\cal O}$, let $\epsilon(\gamma)$ be the sign of $\det(df-1)$, as
in the introduction.
Now define
\begin{equation}
\label{eqn:IEtaDef}
I_\eta \eqdef \prod_{\gamma\in{\cal O}}(1-[\gamma])^{-\epsilon(\gamma)}\det(P)
\in \Lambda.
\end{equation}
By the definition of determinant, the function $I_\eta:H_1(X',\partial
X')\to{\Bbb Z}$ is nonzero only on elements of $H(\eta)$. We now wish to define
$$I\eqdef I_\eta\circ\Psi_\eta:{\cal S}\to{\Bbb Z}.$$
If we change the orientation choices above, this will multiply $I$ by
$\pm1$.

\proposition
\label{pro:invariance}
Modulo the above sign ambiguity, $I$ depends only on the cohomology
class of $\eta$.
\propositionEnd

\proofSketch
Given two cohomologous $\eta$'s, one can show that they are homotopic
through closed forms in the same cohomology class with no index 0 or 3
critical points.  If we generically deform $\eta$ and/or the metric over
time, during a time interval when $\phi$ remains generic, none of the terms
in the definition of $I_\eta$ change.  At certain times, $\eta$ may fail to be
generic in one of the following two ways:
\begin{description}
\item{(a)} Two critical points of index 1 and 2 with a single short flow
line $\gamma$ between them annihilate each other (or are created). 
\item{(b)} There is a flow line $\gamma$ connecting two critical points
with the same index. 
\end{description}

In case (a), if $x$ is an
index 1 point with a flow line to the index 2 point being annihilated, and
if $y$ is an index 2 point with a flow line to the index 1 point being
annihilated, then after the annihilation, these two flow lines fuse into a
single flow line from $x$ to $y$.  In other words, the path matrix $P$ changes
as follows:
$$\left(\begin{array}{cc} \pm[\gamma] & v \\ w & P' \end{array}\right)
\longmapsto
(P'+wv).$$
(Here $v$ and $w$ are the row and the column associated to the index 2 point
and the index 1 point, respectively.)  This divides the determinant, and
hence $I_\eta$, by
$\pm[\gamma]$.  However the isomorphisms $\Psi_{\eta}$ before and after the
annihilation also differ by a factor of $\pm[\gamma]$.  Thus $I$ is changed
only by $\pm1$.

In case (b), suppose without loss that the two critical points $x,y$ have
index one.  Let the flow line $\gamma$ start at $x$ and end at $y$.
Suppose first that $x\neq y$.  The effect of this catastrophe is to replace
$P(x)$ by $P(x)\pm[\gamma] P(y)$.  This does not change $\det(P)$.

When $x=y$, $P(x)$ is multiplied by $(1-[\gamma])^{\pm1}$, but a closed
orbit is created or destroyed to cancel this.  This will not
change $I_\eta$ for the same reason that the left side of
Theorem~\ref{thm:refinement} is invariant, which was remarked upon in the
introduction.

Note that the set of times at which changes of type (b) occur is in general
not discrete.  However, for any integer $k$, the set of times at which
terms $[\gamma]$ in $I_\eta$ with $|\gamma\cdot\Sigma|\le k$ change is
discrete.  So if we discard terms $[\gamma]$ in $I_\eta$ with
$|\gamma\cdot\Sigma|>k$, the resulting expression is invariant.  Taking
$k\to\infty$, it follows that $I_\eta$ is invariant.
\qedNoSkip

\conjecture
\label{cnj:main}
let $X$ be a closed oriented 3-manifold with $b^1(X)>0$. Then $I$ does not
depend on the choice of $\eta\in H^1(X;{\Bbb Z})$, and
$$SW=\pm I.$$
\conjectureEnd

Note that when $b^1(X)=1$,
we define $SW$ to be the limit of the number of
solutions to the equations perturbed by $-ir*\eta$, where $r$ is a large
real number.  (See Meng-Taubes~\cite{meng-taubes}.)

This conjecture is analogous to Taubes' results relating the Seiberg-Witten
and the Gromov invariants in 4 dimensions.  The funny way that closed
orbits are counted in $I$ is analogous to more intricate results in
\cite{taubes:counting}.  The idea of (part of) the proof is that if we have
a nonzero Seiberg-Witten invariant and take the limit as $r\to\infty$, we
get a sequence of Seiberg-Witten solutions such that the zero set of
the $E$ component of the spinor converges to one of the submanifolds
that $I$ counts.

\subsection{The Meng-Taubes formula}
\label{sec:meng-taubes}

\proofOf{Theorem~\ref{thm:application}}
Let $\eta$ be the harmonic 1-form representing $\alpha$, and perturb it
slightly so that $\eta$ is $d$ of a generic Morse function $\phi:X\to S^1$.
Since we started with a harmonic form, there can be no index 0 or 3
critical points.  Let $0$ be a regular value of $\phi$,
and let $\Sigma=\phi^{-1}(0)$.

Let $S$ be a $Spin^c$ structure.  We have
$$\alpha(c_1(\det S))  =  \int_\Sigma c_1(\det S).$$
On $\Sigma$, we have the decomposition
$$S=E\oplus K^{-1}E$$
from \S\ref{sec:conjecture}.  Clearly $K^{-1}|\Sigma=T\Sigma$, so
$$\alpha(c_1(\det S)) = \chi(\Sigma) + 2\Sigma\cdot\Psi_\eta(S).$$
By Conjecture~\ref{cnj:main},
$$\sum_{S\in{\cal S}}SW(S)t^{\Sigma\cdot\Psi_\eta(S)} =
\pm\rho(I_\eta),$$ 
where $\rho:\Lambda \to L_{{\Bbb Z}}$ sends $\gamma\mapsto
t^{\Sigma\cdot\gamma}$.  Thus
$$\sum_{S\in{\cal S}}SW(S)t^{\alpha(c_1(\det S))/2} =
\pm t^{\chi(\Sigma)/2}\rho(I_\eta).$$
To compute $\rho(I_\eta)$, observe that in the notation of \S1,
\begin{eqnarray*}
\rho\left(\prod_{\gamma\in{\cal O}} (1-[\gamma])^{-\epsilon(\gamma)}\right)
& = & \zeta(f),\\
\rho(\det(P)) & = & \det(d:M^1\otimes L_{{\Bbb Q}}\to M^2\otimes
L_{{\Bbb Q}}).
\end{eqnarray*}
Since $M^1$ and $M^2$ are the only nontrivial terms in the Morse complex,
$\det(d:M^1\otimes L_{{\Bbb Q}}\to M^2\otimes L_{{\Bbb Q}})$ is $\tau(M)$ if
$M^*\otimes L_{{\Bbb Q}}$ is acyclic, and zero otherwise.  So by
(\ref{eqn:IEtaDef}), $\rho(I_\eta)$ equals the left side of
Theorem~\ref{thm:refinement} when $M^*\otimes
L_{{\Bbb Q}}$ is acyclic, and zero otherwise.  We are done by
Theorem~\ref{thm:refinement}.
\qedNoSkip

\paragraph{Milnor torsion.}
Let $H=H_1(X)/\mbox{Torsion}=H^2(X;{\Bbb Z})/\mbox{Torsion}$, and let
$\hat{X}$ be the covering of $X$ 
whose monodromy is the projection $\pi_1(X)\to H$.  The ``Milnor torsion''
$MT$ of $X$ is the torsion of the complex $C_*(\hat{X})\otimes
Q({\Bbb Z}[H])$, where $C_*(\hat{X})$ is the cellular complex coming from an
equivariant cell decomposition, and ${\Bbb Z}[H]$ is the group ring of $H$
(where $H$ is written multiplicatively).  This is a well defined element of
$Q({\Bbb Z}[H])/H$, up to sign.  In fact the sign can be specified, by the
same data needed to specify the sign of SW (see  Meng-Taubes
\cite{meng-taubes}, Turaev \cite{turaev}).  The Milnor torsion is
defined to be zero if the complex is not acyclic.

When $b^1(X)>1$, it turns out that $MT\in{\Bbb Z}[H]/H$.
(See Turaev \cite[Thm.\ 1.1.2]{turaev}.)  Furthermore there is a unique
element in this equivalence class invariant under the map that sends $h\to
h^{-1}$ for $h\in H$ \cite[\S1.11.5]{turaev}.  When $b^1(X)>1$ we will
identify $MT$ with this element of ${\Bbb Z}[H]$.

\bigskip

Following Meng-Taubes \cite{meng-taubes}, define
$$\underline{SW}\eqdef \sum_{S\in{\cal S}}SW(S)\frac{c_1(\det
S)}{2}\in{\Bbb Z}[[H]].$$
Here ${\Bbb Z}[[H]]$ is the set of functions $H\to{\Bbb Z}$ that do not
necessarily have finite support.  We can now deduce part of the Meng-Taubes
formula:

\theorem[assuming Conjecture~\ref{cnj:main}]
\label{thm:meng-taubes}
Let $X$ be a closed oriented
3-manifold with $b^1(X)>0$.  Then
$$\underline{SW}=\pm MT.$$
\theoremEnd

\lemma
\label{lem:uniqueness}
Let $G$ be a free abelian group on $m$ generators and let $f,g\in{\Bbb Z}[G]$.
Suppose that  for every homomorphism $\alpha:G\to{\Bbb Z}$,
$\alpha(f)=\alpha(g)$ in ${\Bbb Z}[{\Bbb Z}]$, up to sign.  Then $f=\pm g$.
\lemmaEnd

\proof
Let $\{e_i\}$ be a free basis for $G$.  Choose an integer $N$ such that
$f$ and $g$ are supported in the set $\{\sum a_ie_i\mid |a_i|<N\}$.  Let
$\alpha$ send $e_i$ to $(2N)^i$.  Then for any integer $k$, the hyperplane
$\{x\in G\mid \alpha(x)=k\}$ contains at most one point of union of the
supports of $f$ and $g$.  Apply the hypothesis to this $\alpha$.
\qed

\proofOf{Theorem~\ref{thm:meng-taubes}}
If $b^1(X)=1$ then this is just Theorem~\ref{thm:application}.  Assume
$b^1(X)>1$.  We have already remarked that $MT\in{\Bbb Z}[H]$.  We also have
$\underline{SW}\in{\Bbb Z}[H]$, by the well known a priori bounds on 
the for the Seiberg-Witten equations (see Witten \cite{witten:monopoles}).
(Note that we do not necessarily have $\underline{SW}\in{\Bbb Z}[H]$ when
$b^1=1$, because here we are making a large perturbation to the equations
which destroys the a priori bounds,  and the Seiberg-Witten invariants are
not invariant under perturbation when $b^1(X)=1$.)

If $\alpha\in H^1(X;{\Bbb Z})$ then $\alpha$ extends to a function
${\Bbb Z}[H] \to {\Bbb Z}[{\Bbb Z}] = {\Bbb Z}[t,t^{-1}]$, and the left hand
side of Theorem~\ref{thm:application} is $\alpha(\underline{SW})$, modulo
signs.  On the other hand the right side of Theorem~\ref{thm:application}
is $\alpha(MT)$.  (This is easy when both the complexes involved are
acyclic, and the general case follows from Turaev
\cite[Thm.\ 1.1.3]{turaev} and Corollary~\ref{cor:isomorphism}.)  So
Theorem~\ref{thm:application} says that
$\alpha(\underline{SW})=\alpha(MT)$, modulo signs and powers of $t$.  Since
both $\underline{SW}$ and $MT$ are symmetric,
$\alpha(\underline{SW})=\alpha(MT)$ modulo signs.  We are done by
Lemma~\ref{lem:uniqueness}.
\qedNoSkip

\paragraph{Final remark.}
Under favorable circumstances one can define the torsion of larger
(i.e. not free) abelian coverings of $X$.  (See e.g.\ Fried
\cite{fried:homological}.)  If $X$ is a fibration over $S^1$ (of any
dimension), then  Fried \cite{fried:homological}
shows that the torsion of the universal abelian cover can be identified
with the zeta function, which is our $I$ in the
3-dimensional case.  This is slightly stronger than our result, since a
3-manifold fibered over $S^1$ may have torsion in $H_1$.
We can use an equivariant version of circle-valued Morse theory to extract
more information about the Seiberg-Witten invariants for other 3-manifolds,
and we intend to discuss this in a future paper.

\end{document}